\documentclass[english]{amsart}
\usepackage[T1]{fontenc}
\usepackage[latin1]{inputenc}
\usepackage{babel}

\makeatletter

\providecommand{\LyX}{L\kern-.1667em\lower.25em\hbox{Y}\kern-.125emX\@}

 \theoremstyle{plain}    
 
 \numberwithin{equation}{section} 
 \numberwithin{figure}{section} 
 \newcommand{\lyxaddress}[1]{
   \par {\raggedright #1 
   \vspace{1.4em}
   \noindent\par}
 }

\makeatother
\begin{document}

\newcommand{\fg}{\mathfrak {g}}

\newcommand{\fb}{\mathfrak {b}}

\newcommand{\sD}{\text {\sf D}}

\newcommand{\cB}{\mathcal{B}}

\newcommand{\cF}{\mathcal{F}}

\newcommand{\cL}{\mathcal{L}}

\newcommand{\cN}{\mathcal{N}}

\newcommand{\fu}{\mathfrak {u}}

\newcommand{\A}{\mathbb {A}}

\newcommand{\N}{\mathbb {N}}

\newcommand{\Q}{\mathbb {Q}}

\newcommand{\Z}{\mathbb {Z}}

\newcommand{\Hom}{\mathrm{Hom}}

\newcommand{\Spec}{\mathrm{Spec}\, }

\newcommand{\an}{\mathrm{an}}

\newcommand{\pr}{\mathrm{pr}}

\title{Springer correspondence via p-adic analytic methods}

\maketitle

\subsection*{Introduction}

Despite its alluring title, this note contains hardly anything original:
it is really a footnote to \cite{Ra}, and the author regrets that
his ignorance of representation theory prevented him from including
this material in \emph{loc.cit}. In an effort to make a virtue of
necessity, we hail the lack of originality instead of dissimulating
it: indeed, the main purpose of this article is to show that the usual
proof (via Deligne-Fourier transform) of the Springer correspondence
for semisimple groups over a field \( k \) of positive characteristic,
works \emph{verbatim} even when \( k \) has characteristic zero,
up to replacing Deligne-Fourier transform by the (\( p \)-adic analytic)
Fourier transform of \cite{Ra}, and the usual étale site by the analytic
étale site. We will be rather sketchy, except in the few points where
\( p \)-adic techniques are brought to bear that might be unfamiliar
to an expert of representation theory.

Especially, we will need to assume that \( k \) is complete under
a non-archimedean valuation, but this is not really a restriction. 

Of course, there are better ways (due to Lusztig and others) to obtain
the Springer correspondence in a characteristic-free way; this note
should be regarded as a proof of concept which hopefully might suggest
new applications.

\subsection*{Semisimple groups}

We begin by recalling a few notations and definitions, following \cite[Ch.VI]{KW}.
Let \( k \) be an algebraically closed field of characteristic zero.

We let \( G \) be a semisimple connected algebraic group of dimension
\( n \) over \( \Spec k \), \( W \) the Weyl group of \( G \),
\( B\subset G \) a Borel subgroup, \( U\subset B \) the unipotent
radical of \( B \), \( \fg  \), \( \fb  \) and \( \fu  \) the
Lie algebras of respectively \( G \), \( B \) and \( U \); set
\( b:=\dim (B) \) and \( u:=\dim (U) \). Let moreover \( \cN \subset \fg  \)
be the nilpotent variety, i.e. the Zariski closed subset consisting
of all the nilpotent elements of \( \fg  \). The adjoint action of
\( G \) on \( \fg  \) induces a proper surjective morphism of schemes\[
q_{\cN }:X:=G\times ^{B}\fu \to \cN .\]
It is well-known, and shown in \cite[Ch.VI, \S11, Lemma 1]{KW} that
\( q_{\cN } \) is semi-small; as a consequence, if we let \( m:=\dim \cN  \),
the complex\[
\Psi :=Rq_{\cN *}\overline{\Q }_{\ell ,X}[m](\frac{m}{2})\]
is a \( G \)-equivariant perverse sheaf on the étale site of the
nilpotent variety \( \cN  \) (here \( [m] \) denotes the shift operator
in the derived category, and \( (\frac{m}{2}) \) denotes Tate twist:
notice that \( m \) is always an even number).

On the other hand, the Killing form induces a natural \( G \)-equivariant
isomorphism\[
\fg \simeq \fg ^{*}:=\Hom _{k}(\fg ,k).\]
Under this isomorphism, \( \fb  \) is naturally identified with the
orthogonal \( \fu ^{\bot } \) of the subspace \( \fu  \). We deduce
a map\[
q':Y:=G\times ^{B}\fu ^{\bot }\to \fg ^{*}\]
which is dual to the above \( q_{\cN } \). In particular, the complex\[
\Phi ':=Rq'_{!}\overline{\Q }_{\ell ,Y}[n]\]
is a \( G \)-left-equivariant perverse sheaf on the étale site of
\( \fg ^{*} \), with a \( G \)-action covering the coadjoint action
of \( G \) on \( \fg ^{*} \). In addition, \( \Phi ' \) carries
a natural right \( W \)-action covering the trivial \( W \)-action
on \( \fg ^{*} \)(see \cite[Ch.VI, \S10]{KW}).

\subsection*{Fourier transform}

We recall some results from \cite{Ra}. Here we assume that \( k \)
is complete under a non-archimedean valuation and its residue characteristic
is \( p>0 \), different from \( \ell  \). Let \( S \) be an analytic
variety over \( k \), in the sense of V.Berkovich \cite{Be} or R.Huber
\cite{Hu}; let also \( E \) be an analytic vector bundle over \( S \),
\( E^{*} \) the dual of \( E \), \( p_{1} \) and \( p_{2} \) the
two projections of \( E\times _{S}E^{*} \) onto respectively \( E \)
and \( E^{*} \), and \( \langle \cdot \rangle :E\times _{S}E^{*}\to \A ^{1}_{S} \)
the natural bilinear pairing. Furthermore, let \( \Lambda  \) be
a torsion ring, filtered union of its finite subrings, such that:

\begin{enumerate}
\item \( \ell ^{n}\Lambda =0 \) for sufficiently large \( n \), and 
\item the group of units \( \Lambda ^{\times } \) contains a homomorphic
image of \( \mu _{p^{\infty }} \), the group of all \( p \)-primary
roots of \( 1 \) in \( k \). 
\end{enumerate}
The Fourier transform\[
\cF :\sD ^{b}(E,\Lambda )\to \sD ^{b}(E^{*},\Lambda )\]
is a triangulated functor on the bounded derived category of sheaves
of \( \Lambda  \)-modules on the (analytic) étale site of \( E \),
with values in the bounded derived category of \( \Lambda  \)-modules
on the étale site of \( E^{*} \). As for any transform of its kind,
it is defined as an {}``integral operator''\[
\cF (K^{\bullet }):=Rp_{2!}(p_{1}^{*}K^{\bullet }\otimes \langle \cdot \rangle ^{*}\cL _{\psi })\qquad (\textrm{for all }K^{\bullet }\in \sD ^{b}(E,\Lambda ))\]
whose {}``kernel'' is of the form \( \cL _{\psi }:=\cL \times ^{\mu _{p^{\infty }}}\Lambda  \),
where \( \psi :\mu _{p^{\infty }}\to \Lambda ^{\times } \) is an
injective group homomorphism (that exists by our assumption) and \( \cL  \)
is a certain \( \mu _{p^{^{\infty }}} \)-torsor on the analytic étale
topology of \( \A ^{1}_{S} \) (the sheaf of étale local sections
of the logarithm mapping : see \cite[\S6.1]{Ra} for the full details). 

The main point is the claim that \( \Psi  \) is the Fourier transform
of \( \Phi ' \) (up to a twist), which enables one to transfer the
\( W \)-action from \( \Phi ' \) onto \( \Psi  \). However, before
being able even to state this assertion, we have to take care of a
minor technical nuisance: namely, the formalism of the \( p \)-adic
Fourier transform is (currently) available only for complexes of torsion
modules (and not for complexes of \( \ell  \)-adic sheaves). More
generally, due to some subtle issues, a workable definition of a (triangulated)
derived category of analytic \( \ell  \)-adic sheaves is still missing.
Luckily, for our present purposes, one does not need a full-fledged
formalism of \( \ell  \)-adic complexes; instead, we can proceed
in a more \emph{ad hoc} manner, as follows.

First, let us denote by \( \overline{\Z }_{\ell } \) the integral
closure of \( \Z _{\ell } \) in \( \overline{\Q }_{\ell } \); for
every \( r\in \N  \) the quotient ring \( \Lambda _{r}:=\overline{\Z }_{\ell }/\ell ^{r+1}\overline{\Z }_{\ell } \)
fulfills the conditions (1) and (2) above. 

Recall that the constant \( \ell  \)-adic sheaf \( \overline{\Q }_{\ell } \)
on \( X \) is represented by the inverse system of torsion sheaves
\( \underline{S}:=(\Z /\ell ^{r+1}\Z _{X}\, |\, r\in \N ) \) and
the complex \( \Psi  \) is represented by the corresponding inverse
system \( \underline{T}:=(Rq_{\cN *}\Z /\ell ^{r+1}\Z _{X}[m](\frac{m}{2})\, |\, r\in \N ) \).
Instead of dealing with the finite rings \( \Z /\ell ^{r+1}\Z  \),
we wish to work with coefficient sheaves that are \( \Lambda _{r} \)-modules;
to this aim, it suffices to tensor termwise the inverse system \( \underline{S} \),
which gives the inverse system \( \underline{S'}:=(\Lambda _{r,X}\, |\, r\in \N ) \);
next, instead of applying termwise the functor \( Rq_{\cN *} \) to
the system \( \underline{S} \), we do the same thing with the system
\( \underline{S'} \); this is a rather harmless procedure: indeed,
the functor \( Rq_{\cN *} \)commutes with the change of ring mapping
induced by the (unique) ring homomorphism \( \Z /\ell ^{r+1}\Z \to \Lambda _{r} \),
hence the result is the same as tensoring termwise the inverse system
\( \underline{T} \) by the ring extensions \( \Lambda _{n} \). Likewise,
one replaces \( \Phi ' \) by an inverse system of \( \Lambda _{r} \)-modules.

We are therefore reduced to compare the complexes\[
\Psi _{r}:=Rq_{\cN *}\Lambda _{r,X}\qquad \textrm{and}\qquad \Phi '_{r}:=Rq'_{!}\Lambda _{r,Y}.\]
Now we jump to the analytic étale site of the analytic spaces \( X^{\an } \)
and \( Y^{\an } \), analytification of the \( k \)-schemes \( X \)
and \( Y \); the complexes \( \Psi _{r} \) and \( \Phi '_{r} \)
determine complexes on the étale sites of respectively \( X^{\an } \)
and \( Y^{\an } \), which we shall denote by the same names. Since
the transition to associated analytic étale complexes is compatible
with all cohomological operations, this does not give rise to ambiguities.

\subsection*{Sketch of the proof}

Now the argument proceeds as in \cite[Ch.VI, \S13]{KW} : we let \( \cB :=G/B \);
one has natural closed imbeddings\[
i:X^{\an }\hookrightarrow \fg \times \cB ^{\an }\qquad i':Y^{\an }\hookrightarrow \fg ^{*}\times \cB ^{\an }\]
We regard the analytic space \( \fg \times \cB ^{\an } \) as a trivial
vector bundle over \( \cB ^{\an } \), whose dual bundle is \( \fg ^{*}\times \cB  \)
and notice that \( i \) (resp. \( i' \)) realizes \( X^{\an } \)
(resp. \( Y^{\an } \)) as a sub-bundle of \( \fg \times \cB  \)
(resp. of the dual bundle). We choose a compatible family of characters
\( \psi _{r}:\mu _{p^{\infty }}\to \Lambda _{r}^{\times } \) for
every \( r\in \N  \), and we consider the corresponding Fourier transforms:\[
\cF _{\cB }:\sD ^{b}(\fg ^{*}\times \cB ^{\an },\Lambda _{r})\to \sD ^{b}(\fg \times \cB ^{\an },\Lambda _{r}).\]
Using \cite[Prop.7.2.1]{Ra} we deduce a natural isomorphism:\[
\cF _{\cB }(i'_{*}\Lambda _{r,Y^{\an }}[n](b))\simeq i_{*}\Lambda _{r,X^{\an }}[m].\]
On the other hand, we have also the Fourier transforms on the affine
space \( \fg ^{*} \), regarded as a vector bundle over a point:\[
\cF :\sD ^{b}(\fg ^{*},\Lambda _{r})\to \sD ^{b}(\fg ,\Lambda _{r}).\]
Let \( \pr :\fg \times \cB ^{\an }\to \fg  \) be the projection;
by \cite[Prop.7.1.19]{Ra} we have a natural isomorphism:\[
R\pr _{!}\cF _{\cB }(i'_{*}\Lambda _{r,Y^{\an }}[n])\simeq \cF (Rq'_{!}\Lambda _{r,Y^{\an }}[n]).\]
Let \( i_{\cN }:\cN \to \fg  \) be the closed imbedding; since \( \pr \circ i=i_{\cN }\circ q_{\cN } \)
we finally obtain a compatible system of natural isomorphisms:\[
\cF (\Phi _{r}')\simeq i_{\cN *}\Psi _{r}(-n)\qquad \textrm{for every }r\in \N .\]
From here on it's all downhill : the \( G \)-equivariance of \( \Psi  \)
allows one to decompose it (in the category of perverse sheaves) as
a sum of factors supported on the orbits of the \( G \)-action on
\( \cN  \); on the other hand, the relation just found shows that
the decomposition must also be \( W \)-equivariant and thus we obtain
a map from the set \( \widehat{W} \) of irreducible representations
of \( W \) to the nilpotent conjugacy classes in \( \fg  \), which
turns out to be injective: this is the (Brylinski-)Springer correspondence
(see \cite[Ch.VI, \S\S12-13]{KW} for the complete story).

~\\

\lyxaddress{Lorenzo Ramero\\
Institut de Mathematiques\\
Université de Bordeaux I\\
351, cours de la Liberation\\
33405 Talence cedex\\
FRANCE \\
~\\
ramero@math.u-bordeaux.fr\\
http://www.math.u-bordeaux.fr/\( \sim  \)ramero}
\end{document}